\documentclass[12pt]{amsart}
\usepackage{amsmath,amsthm}
\usepackage{amssymb}
\newtheorem{thm}{Theorem}[section]
\newtheorem*{TA}{Theorem A}
\newtheorem{lemma}[thm]{Lemma}
\newtheorem{prop}[thm]{Proposition}
\newtheorem{cor}[thm]{Corollary}

\theoremstyle{definition}
\newtheorem{defn}[thm]{Definition}
\newtheorem{example}[thm]{Example}

\theoremstyle{remark}
\newtheorem{remark}[thm]{Remark}

\newtheorem{fact}[thm]{Fact}

\begin{document}
%
%
\newcounter{nu}
%
%
\def\sect{\setcounter{thm}{0} \section}
\def\equ{\setcounter{equation}{\value{thm}}\stepcounter{subsection} 
           \begin{equation}}
\def\endeq{\end{equation}\setcounter{thm}{\value{equation}}}
\def\equa{\setcounter{equation}{\value{thm}}\stepcounter{subsection}
           \begin{eqnarray}}
\def\endeqa{\end{eqnarray} \setcounter{thm}{\value{equation}}}
%
%
\newcommand{\ra}{\rightarrow}
\newcommand{\lra}{\longrightarrow}
\newcommand{\epi}{\ra\hspace{-3.2mm}\ra}
\newcommand{\hra}{\hookrightarrow}
\newcommand{\hsm}{\hspace{3 mm}}
\newcommand{\hsp}{\hspace{7 mm}}
\newcommand{\vsm}{\vspace{2 mm}}
\newcommand{\DEF}{:=}
\newcommand{\eps}{\varepsilon}
\newcommand{\rest}[1]{\lvert_{#1}}
\newcommand{\lrc}[1]{\langle\!\langle{#1}\rangle\!\rangle}
\newcommand{\real}[1]{\lvert{#1}\rvert}
\newcommand{\mreal}[1]{\|{#1}\|}
\newcommand{\colim}{\operatorname{colim}}
%
%
\newcommand{\Abgp}{{\mathcal AbGp}}
\newcommand{\CC}{{\mathcal C}}
\newcommand{\del}{\mbox{\scriptsize $\Delta$}}
\newcommand{\FF}{{\mathcal F}}
\newcommand{\GG}{{\mathcal G}}
\newcommand{\Set}{{\mathcal Set}}
\newcommand{\Ss}{{\mathcal S}_{\ast}}
\newcommand{\Sk}{\Ss^{Kan}}
\newcommand{\Sc}{{\mathcal S}_{0}}
\newcommand{\Sr}{{\mathcal S}_{0}^{Kan}}
\newcommand{\Ts}{{\mathcal T}_{\ast}}
\newcommand{\Tc}{{\mathcal T}_{0}}
\newcommand{\WW}[1]{{\mathcal W}_{#1}}
\newcommand{\PAlg}{\hbox{$\Pi$-$Alg$}}
%
%
\newcommand{\N}{\mathbb N}
\newcommand{\R}{\mathbb R}
\newcommand{\ZZ}{\mathbb Z}
\newcommand{\bk}{\mathbf{k}}
\newcommand{\bl}{\mathbf{l}}
\newcommand{\bm}{\mathbf{m}}
\newcommand{\bn}{\mathbf{n}}
\newcommand{\bo}{\mathbf{1}}
%
%
\newcommand{\As}{A_{\ast}}
\newcommand{\Bs}{B_{\ast}}
\newcommand{\Cs}{C_{\ast}}
\newcommand{\Es}{E_{\ast}}
\newcommand{\pis}{\pi_{\ast}}
\newcommand{\Pa}{$\Pi $-algebra}
\newcommand{\Gs}{G_{\ast}}
\newcommand{\Hs}{G'_{\ast}}
%
%
\newcommand{\A}{\mathbf{A}}
\newcommand{\B}{\mathbf{B}}
\newcommand{\C}{\mathbf{C}}
\newcommand{\bD}{{\mathbf \Delta}}
\newcommand{\dD}{\bD_{\partial}}
\newcommand{\Do}{\bD^{op}}
\newcommand{\dDo}{\dD^{op}}
\newcommand{\be}[1]{\mathbf{e}^{#1}}
\newcommand{\G}{\mathbf{G}}
\newcommand{\K}{\mathbf{K}}
\newcommand{\LL}{\mathbf{L}}
\newcommand{\M}{\mathbf{M}}
\newcommand{\map}{\mathbf{map}_{\ast}}
\newcommand{\Pe}[1]{\mathbf{P}^{#1}}
\newcommand{\PP}[2]{\Pe{#1}_{#2}}
\newcommand{\bS}[1]{\mathbf{S}^{#1}}
\newcommand{\U}{\mathbf{U}}
\newcommand{\V}{\mathbf{V}}
\newcommand{\sV}[1]{\V\!_{#1}}
\newcommand{\uV}[2]{\V^{#1}\!\!\!_{#2}\,}
\newcommand{\W}{\mathbf{W}}
\newcommand{\uW}[2]{\W^{#1}_{#2}}
\newcommand{\bW}[2]{\bar{\W}^{#1}\!\!\!_{#2}\,}
\newcommand{\X}{\mathbf{X}}
\newcommand{\Y}{\mathbf{Y}}
\newcommand{\Z}{\mathbf{Z}}
%
%
\newcommand{\pas}{$\Delta$-simplicial}
\newcommand{\pss}{\pas \ space}
\newcommand{\Asd}{A_{\ast\bullet}}
\newcommand{\diag}{\operatorname{diag}}
\newcommand{\Gd}{\G_{\bullet}}
\newcommand{\Ld}{\LL_{\bullet}}
\newcommand{\Md}{\M_{\bullet}}
\newcommand{\Ud}{\U_{\bullet}}
\newcommand{\Vd}{\V_{\bullet}}
\newcommand{\Vud}{\Vd^{\Delta}}
\newcommand{\hVud}{\hat{\V}_{\bullet}^{\Delta}}
\newcommand{\Wd}{\W_{\bullet}}
\newcommand{\hWd}{^{h}\Wd}
\newcommand{\Wud}{\W^{\Delta}_{\bullet}}
\newcommand{\hWud}{^{h}\Wud}
\newcommand{\Xd}{\X_{\bullet}}
\newcommand{\Xud}{\Xd^{\Delta}}
\newcommand{\hXud}{\hat{\X}_{\bullet}^{\Delta}}
\newcommand{\Yd}{\Y_{\bullet}}
%
%
\title{Loop spaces and homotopy operations}
\author{David Blanc} 
\address{University of Haifa, 31905 Haifa, Israel}
\email{blanc@mathcs2.haifa.ac.il}
\date{February 6, 1997}

\begin{abstract}
We describe an obstruction theory for an $H$-space $\X$ to be
a loop space, in terms of higher homotopy operations taking value in \ 
$\pis\X$. \ These depend on first algebraically ``delooping'' the \Pa \ 
$\pis\X$, \ using the $H$-space structure on $\X$, and then trying to realize
the delooped \Pa.
\end{abstract}
\maketitle

%
%
\section{introduction}
\label{cint}
An $H$-space is a topological space $\X$ with a multiplication; 
the motivating example is a topological group $\G$, which from the point of 
view of homotopy theory is just a loop space: \ 
$\G\simeq\Omega B\G=\map(\bS{1},B\G)$.\hsm
The question of whether a given $H$-space $\X$ is, up to homotopy, a
loop space, and thus a topological group (cf.\ \cite[\S 3]{MilnC1}), \
has been studied from a variety of viewpoints \ -- \ see 
\cite{AdSp,BauG,DLaP,FuchsD,HilLS,MayG,StaH1,StaH2,SteM,SugG,ZabH}, and
the surveys in \cite{StaH}, \cite[\S 1]{StaHC}, and \cite[Part II]{KaneH}. 
Here we address this question from the aspect of homotopy operations, in the 
classical sense of operations on homotopy groups.

As is well known, the homotopy groups of a space $\X$ 
have  Whitehead products and composition operations defined on them; 
in addition, there are various higher order operations on \ $\pis\X$, \ 
such as Toda brackets; and the totality of these actually determine the 
homotopy type of $\X$ \ (cf.\ \cite[\S 7.17]{BlaHH}). \ They should thus 
enable us \ -- \ in theory \ -- \ to determine whether $\X$ 
is a loop space, up to homotopy. \ It is the purpose of this note 
to explain in what sense this can actually be done:

First, we show how an $H$-space structure on $\X$ can be used to
define the action of the primary homotopy operations on the shifted
homotopy groups \ $\Gs=\pi_{\ast-1}\X$ \ (which are isomorphic to \ $\pis\Y$ \
if \ $\X\simeq\Omega\Y$). \ This action will behave properly
with respect to composition of operations if $\X$ is homotopy-associative, 
and will lift to a topological action of the monoid of all maps between 
spheres if and only if $\X$ is a loop space (see Theorem \ref{tone} 
below for the precise statement). \ The obstructions 
to having such a topological action may be formulated in the framework of the
obstruction theories for realizing \Pa s and their morphisms described
in \cite{BlaHH}, which are stated in turn in terms of certain higher 
homotopy operations:

\begin{TA}[Theorem \ref{ttwo} below]
An $H$-group $\X$ is $H$-equivalent to a loop space if and only if 
the collection of higher homotopy operations defined in section \ref{crss} 
below (taking value in homotopy groups) vanish coherently.
\end{TA} 

The question of whether a given topological space $\X$ supports an $H$-space 
structure to begin with was addressed in \cite{BlaHO}, where a similar 
obstruction theory, in terms of higher homotopy operations, was defined.

\subsection{notation and conventions}
\label{snac}\stepcounter{thm}

$\Ts$ \ will denote the category of pointed $CW$ complexes
with base-point preserving maps, and by a \textit{space} we shall always
mean an object in $\Ts$, \ which will be denoted by a boldface
letter: \ $\A,\B,\dotsc,\X$, \ $\bS{n}$, \ and so on. The base-point will be 
written \ $\ast\in\X$. \ The full subcategory of $0$-connected spaces
will be denoted by \ $\Tc$. \ 
$\bD[n]$ \ is the standard topological $n$-simplex in \ $\R^{n+1}$.

The space of Moore loops on \ $\Y\in\Tc$ \ will be denoted by \ $\Omega\Y$. \
This is homotopy equivalent to the usual loop space: that is,
the space of pointed maps \ $\map(\bS{1},\Y)$ \ 
(see \cite[III, Corollary 2.19]{GWhE}). \ 
The reduced suspension of $\X$ is denoted by \ $\Sigma\X$.

$\Abgp$ \ is the category of abelian groups, and \ $gr\Abgp$ \ the category of
positively graded abelian groups.

\begin{defn}\label{dso}\stepcounter{subsection}
$\bD$ is the category of ordered sequences \ 
$\bn= \langle 0,1,\ldots,n \rangle$ \ ($n\in \N$), \ with 
order-preserving maps, and \ $\dD$ \ the subcategory having the same 
objects, but allowing only one-to-one morphisms (so in particular, morphisms 
from $\bn$ to $\bm$ exist only for \ $n\leq m$). \ $\Do$, \ $\dDo$ \ are 
the opposite categories.

As usual, a \textit{simplicial object} over any category $\C$  is a functor \ 
$X:\Do\ra\C$; \ more explicitly, it is a sequence of objects \ 
$\{ X_{n} \}_{n=0}^{\infty}$ \ in $\CC$, equipped with \textit{face maps} \ 
$d_{i}:X_{n}\ra X_{n-1}$ \ and \textit{degeneracies} \ 
$s_{j}:X_{n}\ra X_{n+1}$ \ ($0\leq i,j\leq n$), \ satisfying the usual
simplicial identities (\cite[\S 1.1]{MayS}). \ We often denote such a simplical
object by \ $X_{\bullet}$. \ 
The category of simplicial objects over $\CC$ is denoted by \ $s\CC$.

Similarly, a functor \ $X:\dDo\ra\C$ \ is called a \textit{\pas\ object} 
over $\CC$; this is simply a simplicial object without the degeneracies, 
and will usually be written \ $X^{\Delta}_{\bullet}$. \ When \ $\CC=Set$, \ 
these have been variously called $\Delta$-sets, \textit{ss}-sets, or 
\textit{restricted} simplicial sets in the literature (see \cite{RSaD1}). 
The category of \pas\ objects over $\CC$ is denoted by \ $\del\CC$. \ 
We shall usually denote the underlying \pas\ object of a simplicial object \ 
$Y_{\bullet}\in s\CC$ \ by \ $Y^{\Delta}_{\bullet}\in \del\CC$.
\end{defn}

The category of pointed simplicial sets will be denoted by \ $\Ss$ \ 
(rather than \ $s\Set_{\ast}$); \ its objects will be denoted by boldface 
letters \ $\K,\LL,\M\ldots$. \ The subcategory of fibrant simplicial sets 
(Kan complexes) will be denoted by \ $\Sk$, \ and that of reduced Kan 
complexes by \ $\Sr$. \  
$|\K|\in\Ts$ \ will denote the geometric realization of a simplicial
set \ $\K\in\Ss$, \ while \ $S\!\X\in\Sk$ \ will denote the singular
simplicial set associated to a space \ $\X\in\Ts$. \ 
$\GG$ is the category of simplicial groups. \ (See 
\cite[\S\S 3,14,15,17]{MayS} for the definitions). 

For each of the categories \ $\CC=\Ts$, $\Tc$, $\Sk$, $\Sr$, or \ $\GG$, \ 
we will denote by \ $[\X,\Y]_{\CC}$ \ (or simply \ $[\X,\Y]$, \ if there is no
danger of confusion) the set of pointed homotopy classes of maps \ 
$\X\ra\Y$ \ (cf.\ \cite[\S 5]{MayS} and \cite[\S 3]{KanT}). \ 
The constant pointed map will be written \ $c_{\ast}$, \ 
or simply $\ast$.\hsm The homotopy category of $\CC$, whose objects are
those of $\CC$, and whose morphisms are homotopy classes of maps in
$\CC$, will be denoted by \ $ho\CC$. \ The adjoint functors \ $S$ \
and \ $|\cdot|$ \ induce equivalences of categories \ $ho\Ts\approx ho\Sk$; \ 
similarly \ $ho\Sr\approx ho\GG$ \ under the adjoint functors \
$G$, $\bar{W}$ \ (see \S \ref{ssg} below).

\begin{defn}\label{dhs}\stepcounter{subsection}
A \textit{$H$-space structure} for a space \ $\X\in\Ts$ \ is a 
choice of an $H$-\textit{multiplication map} \ $m:\X\times\X\ra\X$ \
such that \ $m\circ i=\nabla$, \  where \ $i:\X\vee\X\hra\X\times\X$ \ 
is the inclusion, and \ $\nabla:\X\vee\X\ra\X$ \ is the fold map 
(induced by the identity on each wedge summand). If $\X$ may be
equipped with such an $m$, we say that \ $\langle \X,m\rangle$ \ 
(or just $\X$) is an \textit{$H$-space}. \ (If we only have \ 
$m\circ i\sim\nabla$, \  we can find a homotopic map \
$m'\sim m$ \ such that \ $m'\circ i=\nabla$, \ since $\X$ is assumed 
to be well-pointed).

An $H$-space \ $\langle \X,m\rangle$ \ is \textit{homotopy-associative} 
if \ $m\circ (m,id_{X})\sim m\circ(id_{X},m):\X\times\X\times\X\ra\X$. \ 
It is an $H$-\textit{group} if it is homotopy-associative and
has a (two-sided) \textit{homotopy inverse} \ $\iota:\X\ra\X$ \ with \
$m\circ (\iota\times id_{X})\circ \Delta\sim c_{\ast} \sim 
m\circ (id_{X}\times\iota)\circ \Delta$, \ 
(where \ $\Delta:\X\ra\X\times\X$ \ is the diagonal). \ In fact, any 
connected homotopy-associative $H$-space is an $H$-group \ 
(cf.\ \cite[X, Theorem 2.2]{GWhE}).

If \ $\langle \X,m\rangle$ \ and \ $\langle \Y,n\rangle$ \ are two
$H$-spaces, a map \ $f:\X\ra\Y$ \ is called an $H$-\textit{map} if \ 
$n\circ(f\times f)\sim f\circ m:\X\times\X\ra\Y$. \ The set of pointed
homotopy classes of $H$-maps \ $\X\ra\Y$ \ will be denoted by \ $[\X,\Y]_{H}$.
\end{defn}

One similarly defines $H$-simplicial sets and simplicial $H$-maps in the 
category \ $\Ss$.

\subsection{organization}
\label{sorg}\stepcounter{thm}

In section \ref{cpa} we review some background material on \Pa s, and 
in section \ref{csps} we explain how an $H$-space structure on $\X$ determines 
the \Pa\ structure of its potential delooping. In section \ref{cspas}
we provide further background on ($\Delta$-)simplicial spaces and \Pa s, and
bisimplicial groups. In section \ref{csg} we show, in the context of 
simplicial groups, that the \Pa\ structure on \ $\pi_{\ast-1}\X$ \ can be made
``topological'' if and only if $\X$ is a loop space (Theorem \ref{tone}). 

Finally, in section \ref{crss} we recall the obstruction theory of 
\cite{BlaHH} for realizing \Pa s in terms of rectifying 
($\Delta$-)simplicial spaces, and explain how it applies to the
recognition of loop spaces (see Theorem \ref{ttwo}). We also simplify the
general obstruction theory in question, by showing that it 
suffices to rectify the underlying \pss\ associated to a free simplicial \Pa\ 
resolution (Proposition \ref{pfour})\vsm .

I would like to thank the referee for his comments.
%
%
\section{\Pa s}
\label{cpa}

In this section we briefly recall some facts on the primary homotopy 
operations.

\begin{defn}\label{dpa}\stepcounter{subsection}
A \textit{\Pa} \ is a graded group \ $\Gs=\{G_{k}\}_{k=1}^{\infty }$ \ 
(abelian in degrees $>1$), together with an action on \ $\Gs$ \ of the 
primary homotopy operations (i.e., compositions and Whitehead products, 
including the ``$\pi_{1}$-action'' of \ $G_{1}$ \ on the higher \ 
$G_{n}$'s, \ as in \cite[X, \S 7]{GWhE}), \ satisfying the usual
universal identities. See \cite[\S 3]{BlaH} or \cite[\S 2.1]{BlaA} 
for a more explicit description. The category of \Pa s (with the obvious 
morphisms) will be denoted \ $\PAlg$.
\end{defn}

\begin{defn}\label{drpa}\stepcounter{subsection}
We say that a space $\X$ \textit{realizes} a \Pa \ $\Gs$ \ 
if there is an isomorphism of \Pa s \ $\Gs\cong \pis\X$. \ (There may 
be non-homotopy equivalent spaces realizing the same \Pa \ -- \ cf.\ 
\cite[\S 7.18]{BlaHH}). Similarly, a morphism of \Pa s \ 
$\phi:\pis \X\ra\pis\Y$ \ (between realizable \Pa s) is \textit{realizable} 
if there is a map \ $f:\X\ra\Y$ \ such that \ $\pis f=\phi$.
\end{defn}

\begin{defn}\label{dfpa}\stepcounter{subsection}
The \textit{free} \Pa s \ are those isomorphic to \ $\pis \W $, \ for 
some (possibly infinite) wedge of spheres $\W$; \ we say that \ 
$\pis\W$ \ is generated by a graded set \ 
$L_{\ast}=\{L_{k}\}_{k=1}^{\infty }$, \ and write \ 
$\pis \W\cong F(L_{\ast})$,  \ if \ 
$\W=\bigvee_{k=1}^{\infty }\bigvee_{{x\in L_{k}}} \bS{k}_{x}$.
\end{defn}

\begin{fact}\label{rfpa}\stepcounter{subsection}
If we let $\Pi $ denote the homotopy category of wedges of spheres, 
and \ $\FF\subset\PAlg$ \ the full subcategory of free \Pa s, then \ 
$\pis:\Pi\ra\FF$ \ is an equivalence of categories. 
Note that any \Pa \ morphism \ $\phi:\Gs\ra\Hs$ \ is uniquely realizable
if \ $\Gs$ \ is a free \Pa.
\end{fact}

For future reference we note the following:

%
%
\begin{lemma}\label{lthree}\stepcounter{subsection}
If \ $\As,\Bs\in\FF$ \ are free \Pa s and \ 
$\As\stackrel{i}{\hookrightarrow} \Bs \stackrel{r}{\epi}\As$ \ is a 
retraction \ ($r\circ i=id_{\As}$), \ then there is a free \Pa \ $\Cs\in\FF$ \ 
such that \ $\Bs=\As\amalg \Cs$.
\end{lemma}

\begin{proof}
Let \ $Q:\PAlg\ra gr\Abgp$ \ be the ``indecomposables'' functor 
(so \ $Q(\pis\W)\cong H_{\ast}(W;\ZZ)$ \ -- \ see \cite[\S 2.2.1]{BlaH}); \ 
then \ $Q(\As)$ \ and \ $Q(\Bs)$ \ are free abelian groups, and since \ 
$Q(\As)\stackrel{Q(i)}{\hookrightarrow} Q(\Bs) \stackrel{Q(r)}{\epi}Q(\As)$ \ 
is a retraction in \ $\Abgp$, \ there is a graded free abelian group \ $\Es$ \ 
such that \ $Q(\Bs)=Q(\As)\oplus \Es$. \ Choosing graded sets of generators \ 
$\{e_{\gamma}\}_{\gamma\in\Gamma}$ \ for \ $\Es$ \ (in degree $1$, choose 
generators for the free group \ $Ker(r)\subseteq \Bs$), \ and setting \ 
$\Cs=F(\{e_{\gamma}\}_{\gamma\in\Gamma})$, \ yields the required decomposition
(by the Hurewicz Theorem).
\end{proof}

\begin{defn}\label{dfpc}\stepcounter{subsection}
Let \ $T:\PAlg\ra\PAlg$ \ be the ``free \Pa'' comonad 
(cf.\ \cite[VI, \S 1]{MacC}), \ defined \ $T\Gs=
\coprod_{k=1}^{\infty}\coprod_{g\in G_{k}\setminus\{0\}}\ \pis\bS{k}_{(g)}$. \ 
The counit \ $\eps=\eps_{\Gs}:T\Gs\epi\Gs$ \ is defined by \ 
$\iota^{k}_{(g)}\mapsto g$ \ (where \ $\iota^{k}_{(g)}$ \ is the canonical 
generator of \ $\pis\bS{k}_{(g)}$), \ and the comultiplication \ 
$\vartheta=\vartheta_{\Gs}:T\Gs\hra T^{2}\Gs$ \ is induced by the 
natural transformation \ $\bar{\vartheta}:id_{\FF}\ra T|_{\FF}$ \ 
defined by \ $x_{k} \mapsto \iota^{k}_{(x_{k})}$.
\end{defn}

\begin{defn}\label{dapa}\stepcounter{subsection}
An \textit{abelian} \Pa \ is one for which all Whitehead products vanish.
\end{defn}

These are indeed the abelian objects of \ $\PAlg$ \ -- \ see \cite[\S 2]{BlaA}.
If $\X$ is an $H$-space, then \ $\pis\X$ \ is an abelian \Pa \ 
(cf.\ \cite[X, (7.8)]{GWhE}).

%
%
\section{Secondary \Pa \ structure}
\label{csps}

We now describe how an $H$-space structure on $\X$ determines the \Pa\
structure of a (potential) classifying space.

\subsection{the James construction}
\label{sjc}\stepcounter{thm}

For any \ $\X\in\Ts$, \ let \ $J\X$ \ be the James reduced product 
construction, with \ $\lambda:J\X\to \Omega\Sigma\X$ \ the 
homotopy equivalence of \cite[VII, (2.6)]{GWhE}, \ and \ 
$j_{X}:\X\hra J\X$ \ and \ $i_{X}:\X \hra\Omega\Sigma\X$ \ 
the natural inclusions.

If \ $\langle \X,m\rangle$ \ is an $H$-space, \ there is a retraction 
of spaces \ $\bar{m}:J\X\ra\X$ \ (with \ $\bar{m}\circ j_{X}=id_{X}$), \ 
defined \ 
\setcounter{equation}{\value{thm}}\stepcounter{subsection}
\begin{equation}\label{eone}
\bar{m}(x_{1},x_{2},\ldots,x_{n}) = 
m(\ldots m(m(x_{1},x_{2}),x_{3}),\ldots, x_{n}) 
\end{equation}
\setcounter{thm}{\value{equation}}

\noindent (cf.\ \cite[Theorem 1.8]{JamR}).

\begin{defn}\label{ddc}\stepcounter{subsection}
Let $\X$ be an $H$-space. \ Given homotopy classes of maps \ 
$\alpha\in[\Sigma\A,\Sigma\B]$ \ and \ $\beta\in[\B,\X]$, \ 
we define the \textit{derived composition} \ 
$\alpha\star\beta \in [\A,\X]$ \ as follows:

Choose representatives \ $f:\Sigma\A\ra\Sigma \B$ \ and \ $g:\B\ra\X$ \ 
for $\alpha,\beta$ \ respectively, and let \
$\lambda^{-1}:\Omega\Sigma\X\to J\X$ \ be any homotopy inverse 
to $\lambda$. \ Then \ $\alpha\star\beta $ \ is represented by the
composite
$$
\A \stackrel{i_{A}}{\lra} \Omega\Sigma\A \stackrel{\Omega\alpha}{\lra}
\Omega\Sigma\B \stackrel{\lambda^{-1}}{\lra} J\B
\stackrel{J\beta}{\lra} J\X \stackrel{\bar{m}}{\lra} \X.
$$
\end{defn}

\begin{fact}\label{fse}\stepcounter{subsection}
Note that if \ $\alpha=\Sigma\bar{\alpha}$ \ for some \
$\bar{\alpha}:\A\ra\B$, \ then \
$\alpha\star\beta=\bar{\alpha}^{\#}\beta$ \ (this is well-defined,
because $\X$ is an $H$-space)\vsm .
\end{fact}

We shall be interested in the case where $\B$ is a wedge of spheres
and \ $\A=\bS{n}$, \ so $\star$ associates a class \ 
$\omega\star(\beta_{1},\ldots,\beta_{k})\in\pi_{n}\X$ \ to any 
$k$-ary homotopy operation \ 
$\omega^{\#}:\pi_{n_{1}+1}(-)\times\ldots\times\pi_{n_{k}+1}(-)
\to\pi_{n+1}(-)$ \ 
and collection of elements \ $\beta_{i}\in\pi_{n_{i}}\X$ \ ($i=1,\ldots,k$).

In particular, if \ $\omega:\bS{p+q+1}\ra\bS{p+1}\vee\bS{q+1}$ \
represents the Whitehead product, one may define a ``Samelson
product'' \ $\omega\star(-,-):\pi_{p}\X \times \pi_{q}\X\ra\pi_{p+q}\X$ \ 
for any $H$-space $\X$, even without assuming associativity or the
existence of a homotopy inverse (compare \cite[X, \S 5]{GWhE}).

However, in general this \ $\omega\star(-,-)$ \ need not enjoy any 
of the usual properties of the Samelson product (bi-additivity,
graded-commutativity, Jacobi identity \ -- \ cf.\ 
\cite[X, Theorems 5.1 \& 5.4]{GWhE}). \ To ensure that they hold, one 
needs further assumptions on $\X$\vsm .

First, we note the following homotopy version of \cite[VII, Theorem 2.5]{GWhE},
which appears to be folklore:
 
%
%
\begin{lemma}\label{lone}\stepcounter{subsection}
If \ $\langle \X,m\rangle$ \ is a homotopy-associative $H$-space, then
any map \ $f:\A\ra\X$ \ extends to an $H$-map \ $\hat{f}:J\A\ra\X$, \ 
which is unique up to homotopy.
\end{lemma}

\begin{proof}
Given \ $f:\A\ra\X$, \ define \ $\hat{f}:J\A\ra\X$ \ by \ 
$$
\hat{f} (x_{1},\ldots,x_{r}) =
m(\ldots m(m(f(x_{1}),f(x_{2})),f(x_{3})),\ldots, f(x_{r})).
$$

\noindent This is an $H$-map by \cite[Lemma 1.4]{NeiP}. \ Now let \
$\hat{g}:J\A\ra\X$ \ be another $H$-map, with a homotopy \ 
$H:f\simeq g\DEF\hat{g}\circ j_{A}$. \ Since $\hat{g}$ is an $H$-map, there 
is a homotopy \ $G:n\circ (\hat{g}\times\hat{g})\simeq \hat{g}\circ m$ \ 
(where \ $m:J\A\times J\A\ra J\A$ \ is the $H$-multiplication). \
Moreover, by \cite[Lemma 1.3(a)]{NeiP} we may assume $G$ is
stationary on \ $J\A\vee J\A$.

For each \ $r\geq 0$, \ let \ $J_{r}\A$ \ denote the $r$-th stage \ 
in the construction of \ $J\A$, \ with \ $j^{s}_{r}:J_{s}\A\hra J_{r}\A$ \ 
and \ $j^{s}:J_{s}\A\hra J\A$ \ the inclusions, starting with \ 
$J_{0}\A=\ast$ \ and \ $J_{1}\A=\A$. \ 
We define \ $T_{r}\A$ \ to be the pushout in the following diagram:

%
%
\begin{picture}(100,90)(-160,-15)

\put(5,50){$J_{r-1}\A$}
\put(40,55){\vector(1,0){28}}
\put(48,60){$i_{1}$}
\put(70,50){$J_{r-1}\A\times\A$}
\put(20,46){\vector(0,-1){33}}
\put(-3,28){$j^{r-1}_{r}$}
\put(105,46){\vector(0,-1){33}}
\put(108,28){$\bar{q}_{r}$}
\put(10,0){$J_{r}\A$}
\put(33,5){\vector(1,0){60}}
\put(60,-5){$\iota$}
\put(95,0){$T_{r}\A$}
\put(75,12){\framebox(15,10){\scriptsize PO}}
\end{picture}

\noindent for \ $r\geq 1$ \ (so \ $T_{1}\A=\A\vee\A$); \
then \ $J_{r+1}\A$ \ is the pushout in:

%
%
\begin{picture}(200,100)(-100,-15)

\put(30,50){$T_{r}\A$}
\put(55,55){\vector(1,0){113}}
\put(65,62){$\psi_{r}=(i_{1},j^{r-1}_{r}\times id)$}
\put(170,50){$J_{r}\A\times\A$}
\put(42,46){\vector(0,-1){33}}
\put(-25,28){$\varphi_{r}=(id,\bar{q}_{r})$}
\put(195,46){\vector(0,-1){33}}
\put(198,28){$\bar{q}_{r+1}$}
\put(30,0){$J_{r}\A$}
\put(55,5){\vector(1,0){123}}
\put(100,-6){$j^{r}_{r+1}$}
\put(180,0){$J_{r+1}\A$}
\put(160,15){\framebox(15,10){\scriptsize PO}}
\end{picture}

Now let \ $\hat{f}_{r}=\hat{f}|_{J_{r}\A}$ \ and \ 
$\hat{g}_{r}=\hat{g}|_{J_{r}\A}$; \ we shall extend \ $H:f\simeq g$ \ 
to a homotopy \ $\hat{H}:\hat{f}\simeq\hat{g}$ \ by inductively
constructing homotopies \ $\hat{H}_{r}:\hat{f}_{r}\simeq\hat{g}_{r}$ \
(starting with \ $\hat{H}_{1}=H$) \ such that \
$\hat{H}_{r}|_{J_{r-1}\A}=\hat{H}_{r-1}$:  \ 
let \ $n_{r}:\X^{r}\ra\X$ \ denote the $n$-fold multiplication \ 
$n_{r}(x_{1},\ldots,x_{r})=n(\ldots n(x_{1},x_{2}),\ldots),x_{r})$ \ 
and \ $q_{r}:\A^{r}\ra J_{r}\A$ \ the quotient map, so that \ 
$n^{r}\circ f^{r}=\hat{f}_{r}\circ q_{r}$.

As a first approximation, define \ 
$\bar{H}_{r+1}:\hat{f}_{r}\times f\simeq\hat{g}_{r}\times g$ \ on \
$J_{r}\A\times\A$ \ in the above pushout to be the sum of homotopies \ 
$\bar{H}_{r+1} = n \circ (\hat{H}_{r}\times H) + G\circ(j^{r}\times j_{A})$. \ 
This does not quite agree with \ $\hat{H}_{r}\circ\varphi_{r}$ \ on \
$T_{r}\A$, \ but since $G$ is stationary on \ $J\A\vee J\A$ \ we have \ 
$\bar{H}_{r+1}|_{J_{r}\A}=n\circ(\hat{H}_{r}\times id) +\mbox{(stationary)}= 
\hat{H}_{r}+\mbox{(stationary)}$ \ and \ 
$\bar{H}_{r+1}|_{J_{r-1}\A\times\A}=n\circ(\hat{H}_{r-1}\times H) +
G\circ (j^{r-1}\times j_{A})=\bar{H}_{r}$. \ 

Since \ $\hat{H}_{1}=H$, \ we see that \ 
$\bar{H}_{2}|_{T_{2}\A}=(H+\mbox{(stationary)},H)$, \ while \ 
$H\circ\varphi_{1} = (H,H)$. \ Thus we may assume by induction that
there is a homotopy of homotopies \ 
$F:\bar{H}_{r+1}|_{T_{r}\A}\simeq \hat{H}_{r}\circ \varphi_{r}$. \ 
Since \ $T_{r}\A\hra J_{r}\A\times\A$ \ is a cofibration, the inclusion
$$
T_{r}\A\times I^{2}\cup(J_{r}\A\times\A)\times(\{0,1\}\times I\cup
I\times\{0\}) \hra (J_{r}\A\times\A)\times I^{2}
$$

\noindent is a trivial cofibration, and thus we may use the homotopy
extension property to obtain a new homotopy \ $\tilde{F}$ \ on \ 
$(J_{r}\A\times\A)\times I^{2}$ \ which restricts to \ 
$\tilde{H}_{r+1}:\hat{f}_{r}\times f\simeq \hat{g}_{r}\times g$ \ 
on \ $J_{r}\A\times\A\times I\times\{1\}$, \ such that \ 
$\tilde{H}_{r+1}$ \ extends \ $\hat{H}_{r}\circ\varphi_{r}$, \ and
thus may be combined with \ $\hat{H}_{r}$ \ to define a homotopy \ 
$\hat{H}_{r+1}$ \ as required.
\end{proof}
%
%
\begin{cor}\label{ctwo}\stepcounter{subsection}
If $\X$ is a homotopy-associative $H$-space, then for any \ $\A\in\Ts$ \ 
the inclusion \ $j_{A}:\A\ra J\A$ \ induces a bijection \ 
$j_{A}^{\ast}:[J\A,X]_{H}\stackrel{\cong}{\lra} [\A,\X]_{\Ts}$.
\end{cor}

\begin{proof}
Since $\X$ is homotopy-associative $H$-space, the retraction \ 
$\bar{n}=\widehat{id}_{X}:J\X\ra\X$ \ is an $H$-map, by the Lemma, so 
we may define \ $\phi:[\A,\X]_{\Ts}\ra[J\A,\X]_{H}$ \ by \ 
$\phi([f])=[\bar{m}\circ J(f)]$, \ and clearly \ 
$j_{A}^{\ast}(\phi([f])=[\bar{m}\circ J(f)\circ j_{A}]=[f]$.\hsm
On the other hand, given an $H$-map \ $g:J\A\to \X$ \ we have \ 
$\bar{m}\circ J(g\circ j_{\A})\circ j_{A}\simeq g\circ j_{A}$, \ which
implies that \ $\bar{m}\circ J(g\circ j_{\A})\simeq g$ \ by the Lemma 
again. \ Thus also \ $\phi(j_{A}^{\ast}([g]))=[g]$.
\end{proof}

\subsection{notation}
\label{npi}\stepcounter{thm}
If $\X$ is a homotopy-associative $H$-space, we shall write \
$\pi_{t}^{H}\X$ \ for \ 
$[\Omega\bS{t},\X]_{H}=[J\bS{t-1},\X]_{H}\cong \pi_{t-1}\X$. 

%
%
\begin{prop}\label{pone}\stepcounter{subsection}
If $\X$ is a homotopy-associative $H$-space, then \ 
$\alpha\star(\beta\star\gamma)=(\alpha^{\#}\beta)\star\gamma$ \ for
any \ $\alpha\in[\Sigma\A,\Sigma\B]$, \ $\beta\in[\Sigma\B,\Sigma\C]$, \ 
and \ $\gamma\in[\C,\X]$.
\end{prop}

\begin{proof}
It suffices to consider \ $\alpha=id_{\Sigma B}$, \ and so to show that

%
%
\begin{picture}(150,90)(-130,-10)
%
%
\put(0,45){$\Omega\Sigma\B$}
\put(30,50){\vector(1,0){63}}
\put(45,53){$\Omega\beta$}
\put(95,45){$\Omega\Sigma\C$}
\put(125,50){\vector(1,0){38}}
\put(135,55){$\Omega\Sigma\gamma$}
\put(165,45){$\Omega\Sigma\X$}
\put(180,41){\vector(0,-1){23}}
\put(182,25){$\hat{m}$}
%
%
\put(30,40){\vector(3,-1){75}}
\put(8,15){$\Omega\Sigma(\beta\star\gamma)$}
\put(110,5){$\Omega\Sigma \X$}
\put(140,10){\vector(1,0){33}}
\put(150,-2){$\hat{m}$}
\put(175,5){$\X$}
\end{picture}

\noindent commutes up to homotopy (where $\hat{m}$ is the composite \ 
$\Omega\Sigma\X\stackrel{\lambda^{-1}}{\lra}J\X\stackrel{\bar{m}}{\lra}\X$) \ 
 \ --  \ or, since \ $\beta\star\gamma$ \ is defined to be the composite \ 
$\hat{m}\circ \Omega\Sigma\gamma \circ \Omega\beta \circ i_{B}$, \ 
that the two composites \
$\phi=\hat{m}\circ\Omega\Sigma\gamma\circ\Omega\beta$ \ and \ 
$\psi=\hat{m}\circ\Omega\Sigma\hat{m}\circ(\Omega\Sigma)^{2}\gamma 
\circ \Omega\Sigma\Omega\beta \circ \Omega\Sigma i_{B}$ \ are homotopic.

Now if $\X$ is a homotopy-associative $H$-space, then $\hat{m}$ is an 
$H$-map by Lemma \ref{lone}, so \ $\phi,\psi:\Omega\Sigma\B\ra\X$ \ 
are $H$-maps. \ By Corollary \ref{ctwo}
it suffices to check that \ $\phi\circ i_{B}\sim \psi\circ i_{B}$ \ -- \
i.e., that \ $\hat{m}\circ\Omega\Sigma\gamma\circ\Omega\beta\circ i_{B}$ \ is
homotopic to the composition of
$$
\B \stackrel{i_{B}}{\lra} \Omega\Sigma\B
\stackrel{\Omega\Sigma i_{B}}{\lra} (\Omega\Sigma)^{2}\B
\stackrel{\Omega\Sigma\Omega\beta}{\lra} (\Omega\Sigma)^{2}\C
\stackrel{(\Omega\Sigma)^{2}\gamma}{\lra} (\Omega\Sigma)^{2}\X
\stackrel{\Omega\Sigma\hat{m}}{\lra} \Omega\Sigma\X
\stackrel{\hat{m}}{\lra} \X.
$$

But \ $\Omega\Sigma\gamma\circ\Omega\beta\circ i_{B}$ \ is adjoint to \ 
$(\Sigma \gamma)\circ\beta$, \ while the composition of
$$
\B \stackrel{i_{B}}{\lra} \Omega\Sigma\B
\stackrel{\Omega\Sigma i_{B}}{\lra} (\Omega\Sigma)^{2}\B
\stackrel{\Omega\Sigma\Omega\beta}{\lra} (\Omega\Sigma)^{2}\C
\stackrel{(\Omega\Sigma)^{2}\gamma}{\lra} (\Omega\Sigma)^{2}\X
\stackrel{\Omega\Sigma\hat{m}}{\lra} \Omega\Sigma\X
$$

\noindent is adjoint to \ 
$\Sigma(\hat{m}\circ\Omega\Sigma\gamma \circ\Omega\beta\circ i_{B})$ \ which
is equal to \ $\Sigma(\hat{m}\circ\widetilde{(\Sigma\gamma) \circ\beta})$, \
(where $\tilde{f}$ denotes the adjoint of $f$). Since for any \ 
$f:\Y\to\Z$ \ the adjoint of \ $\Sigma f$ \ is \ 
$\Omega\Sigma f\circ i_{Y}$, \ we see \ 
$\hat{m}\circ \widetilde{\Sigma f}\sim f$, \ which completes the proof.
\end{proof}

It is readily verified that when \ $\X\simeq\Omega\Y$, \ the secondary
composition is the adjoint of the usual composition in \ $\pis\Y$; \
thus we have:
%
%
\begin{cor}\label{cone}\stepcounter{subsection}
If $\X$ is an $H$-group, then the graded abelian group \ $\Gs$, \ 
defined by \ $G_{k}=\pi^{H}_{k}\X\cong \pi_{k-1}\X$ \ (with \
$\bar{\gamma}\in G_{k}$ \ corresponding to \ $\gamma\in\pi_{k-1}\X$), \ 
has a \Pa\ structure defined by the derived compositions:\hsm 
that is, \ if \ $\psi\in\pi_{k}(\bS{t_{1}}\vee\ldots\vee\bS{t_{n}})$ \
and \ $\bar{\gamma}_{j}\in G_{t_{j}}$ \ for \ $1\leq j\leq n$, \ 
then \ 
$$
\psi^{\#}(\bar{\gamma}_{1},\ldots, \bar{\gamma}_{n})\DEF 
\overline{\psi\star(\gamma_{1},\ldots, \gamma_{n})}\in G_{k}.
$$
If \ $\X\simeq\Omega\Y$, \ then \ $\Gs$ \ is isomorphic as a \Pa\ to \ 
$\pis\Y$.
\end{cor}

\begin{defn}\label{ddlp}\stepcounter{subsection}
For any $H$-group \ $\langle\X,m\rangle$, \ the \Pa\ structure on the 
graded abelian group\ $\Gs$ \ of Corollary \ref{cone} will be called the 
\textit{delooping} of \ $\pis\X$, \ and denoted by \ 
$\Omega^{-1}\pis\X$ \ (so in particular \ 
$\Omega^{-1}\pis\Omega\Y\cong\pis\Y$).
\end{defn}

\begin{remark}\label{rdlp}\stepcounter{subsection}
Note that Corollary \ref{cone} provides us with an algebraic
obstruction to delooping a space $\X$: \ if there is no way of putting
a \Pa\ structure on the graded abelian group \ $\Gs=\pi_{\ast-1}\X$ \
which is consistent with Fact \ref{fse}, then $\X$ is not a loop
space, or even a homotopy-associative $H$-space. (This is of course
assuming that the \Pa\ \ $\pis\X$ \ is abelian \ -- \ otherwise $\X$
cannot even be an $H$-space.)
\end{remark}

\begin{example}\label{edlp}\stepcounter{subsection}
Consider the \Pa\ \ $\Gs$ \ defined by \ $G_{2}=\ZZ\langle x\rangle$, \ 
(i.e., $x$ generates the cyclic group \ $G_{2}$), \ 
$G_{3}=\ZZ/2\langle \eta_{2}^{\#}x\rangle$, \ 
$G_{4}=\ZZ/2\langle \eta_{3}^{\#}\eta_{2}^{\#} x\rangle$, \ and \ 
$G_{5}=\ZZ/2\langle\eta_{4}^{\#}\eta_{3}^{\#}\eta_{2}^{\#} x\rangle$, \ 
with \ $G_{t}=0$ \ for \ $t\neq 2,3,4,5$ \ and all Whitehead products zero.

There can be no homotopy-associative $H$-space $\X$ with \
$\pis\X\cong\Gs$, \ since the \Pa\ \ $\Gs'=\Omega^{-1}\Gs$ \ cannot be
defined consistently: \ we would have \ 
$G'_{3}=\ZZ\langle\bar{x}\rangle$, \ 
$G'_{4}=\ZZ/2\langle \eta_{3}^{\#}\bar{x}\rangle$, \ 
$G'_{5}=\ZZ/2\langle \eta_{4}^{\#}\eta_{3}^{\#} \bar{x}\rangle$, \ and \ 
$G'_{6}=\ZZ/2\langle\eta_{5}^{\#}\eta_{4}^{\#}\eta_{3}^{\#}\bar{x}\rangle$ \ 
by Fact \ref{fse}; \ but \ 
$\pi_{6}\bS{3}=\ZZ/12\langle\alpha\rangle$ \ with \ 
$6\alpha=\eta_{5}^{\#}\eta_{4}^{\#}\eta_{3}$, \ and thus \
$\alpha^{\#}\bar{x}\in G'_{6}$ \ cannot be defined consistently with
the fact that \ $(6\alpha)^{\#}\bar{x}\neq 0$\vsm .

(We do not claim that \ $\Gs$ \ is realizable; but the obstructions to
realizing \ $\Gs$ \  by a space \ $\X\in\Ts$ \ require secondary 
(or higher order) information, while the obstructions to its
realization by an $H$-group are primary.)
\end{example}

%
%
\section{Simplicial spaces and \Pa s}
\label{cspas}

We next recall some background on ($\Delta$-)simplicial spaces and \Pa s, and
bisimplicial groups:

\begin{defn}\label{daso}\stepcounter{subsection}
Recall that a simplicial object over a category $\C$  is a functor \ 
$X:\Do\ra\C$ \ (\S \ref{dso}); an \textit{augmented} simplical object \ 
$X_{\bullet}\ra A$ \ over $\CC$ is a simplicial object \ 
$X_{\bullet}\in s\CC$, \ together with an augmentation \
$\varepsilon:X_{0}\ra A$ \ in $\CC$ such that \ 
%
\setcounter{equation}{\value{thm}}\stepcounter{subsection} 
\begin{equation}\label{eeight}
\varepsilon\circ d_{1}=\varepsilon\circ d_{0}. \ 
\end{equation}
\setcounter{thm}{\value{equation}}
\noindent Similarly for an \textit{augmented} \pas\ object. 
\end{defn}

\begin{defn}\label{dfsr}\stepcounter{subsection}
A simplicial \Pa \ $\Asd$ \ is called \textit{free} if for each \
$n\geq 0$ \ there is a graded set \ $T^{n}\subseteq (\As)_{n}$ \
such that \ $\As)_{n}$ \ is the free \Pa \ generated by \ $T^{n}$, \
and each degeneracy map \ $s_{j} :(\As)_{n}\ra (\As)_{n+1}$ \
takes \ $T^{n}$ \ to \ $T^{n+1}$.

A \textit{free simplicial resolution} of a \Pa \ $\Gs$ \ 
is defined to be an augmented simplicial \Pa  \ $\Asd\ra\Gs $, \ such
that \ $\Asd$ \ is a free simplicial \Pa, \ the homotopy groups of the 
simplicial group \ $A_{k\bullet}$ \ vanish in dimensions \ $n\geq 1$, \ 
and the augmentation induces an isomorphism \ 
$\pi _{0}A_{k\bullet}\cong G_{k}$.
\end{defn}

Such resolutions always exist, for any \Pa \ $\Gs$ \ -- \ see
\cite[II, \S 4]{QuH}, or the explicit construction in \cite[\S 4.3]{BlaH}.

\subsection{realization}
\label{sreal}\stepcounter{thm}

Let \ $\Wd\in s\Ts$ \ be a simplicial space: \ its \textit{realization} 
(or homotopy colimit) is a space \ $\real{\Wd}\in\Ts$ \ constructed by 
making identifications in \ $\coprod_{n=0}^{\infty}\W_{n}\times\bD[n]$ \
according to the face and degeneracy maps of \ $\Wd$ \ 
(cf.\ \cite[\S 1]{SegC}). \ There is also a \textit{modified realization} \ 
$\mreal{\Wd}\in\Ts$, \ defined similarly, but without making the 
identifications along the degeneracies;  for ``good'' simplicial spaces 
(which include all those we shall consider here) one has \ 
$\mreal{\Wd}\xrightarrow{\simeq}\real{\Wd}$ \ (cf.\ \cite[App.\ A]{SegCC}).
Of course, \ $\mreal{\Wud}$ \ is also defined for \pss s \ $\Wud\in \del\Ts$.

For any reasonable simplicial space \ $\Wd$, \ there is a first 
quadrant spectral sequence with
\setcounter{equation}{\value{thm}}\stepcounter{subsection} 
\begin{equation}\label{ethree}
E^{2}_{s,t}=\pi_{s}(\pi_{t}\Wd) \Rightarrow \pi_{s+t}\real{\Wd}
\end{equation}
\setcounter{thm}{\value{equation}}

\noindent (see \cite[Thm B.5]{BFrH} and \cite[App.]{BLodV}).

\begin{defn}\label{drxs}\stepcounter{subsection}
For any connected \ $\X\in\Ts$, \ an augmented simplical space \ 
$\Wd\ra\X$ \ is called a \textit{resolution of $\X$ by spheres\/} 
if each $\W_{n}$ \ is homotopy equivalent to a wedge of spheres, 
and \ $\pis\Wd\ra\pis\X$ \ is a  free simplicial resolution of \Pa s 
(Def.\ \ref{dfsr}).

Using the above spectral sequence, we see that the natural map \ 
$\W_{0}\ra \real{\Wd}$ \ then induces an isomorphism \ 
$\pis\X\cong\pis\real{\Wd}$, \ so \ $\real{\Wd}\simeq\X$.
\end{defn}

%
%
\section{A simplicial group version}
\label{csg}

For our purposes it will be convenient to work at times in the
category $\GG$ of simplicial groups.   First, we recall some basic
definitions and facts:

\subsection{simplicial groups}
\label{ssg}\stepcounter{thm}

Let \ $F:\Ss\ra\GG$ \ denote the free group functor of \cite[\S
2]{MilnF}; \ this is the simplicial version of the James construction,
and in particular \ $|F\K|\simeq J|\K|$. \ 

Let \ $G:\Ss\ra\GG$ \ be Kan's simplicial loop functor \ (cf.\ 
\cite[Def.\ 26.3]{MayS}), \ with \ $\bar{W}:\GG\ra\Sr$ \ its adjoint,
the Eilenberg-Mac Lane classifying space functor 
(cf.\ \cite[\S 21]{MayS}). \ 

Then \ $|G\K|\simeq\Omega|\K|$ \ and \ $|\K|\simeq |\bar{W}G\K|$. \ 
Moreover, unlike \ $\Ts$, \ where we have only a (weak) homotopy 
equivalence, in $\GG$ there is a canonical isomorphism \ 
$\phi:F\K\cong G\Sigma\K$ \ (cf.\ \cite[Prop.\ 4.15]{CurtS}), 
and there are natural bijections
\setcounter{equation}{\value{thm}}\stepcounter{subsection}
\begin{equation}\label{etwo}
Hom_{\Ss}(\Sigma\LL, \bar{W}F\K)~\cong\ 
Hom_{\GG}(G\Sigma\LL, F\K)~\stackrel{\phi^{\ast}}{\to}\ 
Hom_{\GG}(F\LL, F\K)~\cong\ 
Hom_{\Ss}(\LL, F\K)
\end{equation}
\setcounter{thm}{\value{equation}}
for any \ $\LL\in\Ss$ \ (induced by the adjunctions), and similarly for 
homotopy classes of maps. 

Thus, we may think of \ $F\bS{n}$ \ as the simplicial group analogue
of the $n$-sphere; \ in particular, if $\K$ is in $\GG$, \ or even if 
$\K$ is just an associative 
$H$-simplicial set which is a Kan complex, we shall write \
$\pi_{t}^{H}\K$ \ for \ $[F\bS{t-1},\K]_{H}$ \ (compare \S
\ref{npi}). \ Similarly, \ $F\be{n}$ \ is the $\GG$ analogue of the 
$n$-disc in the sense that any nullhomotopic map \ $f:F\bS{n-1}\ra\K$ \  
extends to \ $F\be{n}$.

\begin{remark}\label{rsg}\stepcounter{subsection}
The same facts as in \S \ref{sreal} hold also if we consider
bisimplicial groups (which we shall think of as simplicial objects \ 
$\Gd\in s\GG$) \ instead of simplicial spaces. \ In this case the 
realization \ $\real{\Wd}$ \ should be replaced by the diagonal \ 
$\diag(\Gd)$, \ and the spectral sequence corresponding to 
\eqref{ethree}, with
\setcounter{equation}{\value{thm}}\stepcounter{subsection} 
\begin{equation}\label{efour}
E^{2}_{s,t}=\pi_{s}(\pi_{t}\Gd) \Rightarrow \pi_{s+t}\diag{\Gd},
\end{equation}
\setcounter{thm}{\value{equation}}
is due to Quillen (cf.\ \cite{QuS})\vsm .
\end{remark}

The above definitions provide us with a functorial simplicial version of the 
derived composition of \S \ref{ddc}:

\begin{defn}\label{dsv}\stepcounter{subsection}
If \ $\K\in\Sk$  \ is an $H$-simplicial set which is a Kan
complex, one again has a retraction  of simplicial sets \ 
$\bar{m}:F\K\ra\K$, \ defined as in \ \eqref{eone}.\hsm
Given a homomorphism of simplicial groups \ $f:F\A\ra F\B$ \ 
and a map of simplicial sets \ $g:\B\ra\K$, \ the composite \ 
$\bar{m}\circ Fg\circ f:F\A\ra\K$  \ will be denoted by \ $f\star g$. 

Note that if \ $\tilde{f}:\Sigma\A\ra\bar{W}F\B$ \ and \ 
$\bar{f}:\A\ra F\B$ \ correspond to $f$ under \eqref{etwo}, the composite \ 
$\bar{m}\circ Fg\circ \bar{f}$ \ corresponds to \ $f\star g$, \ 
and represents the derived composition \ 
$[\bar{f}]\star[g]$ \ in \ $[\A,\K]_{\Ss}\cong [|\A|,|\K|]_{\Ts}$.
\end{defn}

\begin{remark}\label{rsv}\stepcounter{subsection}
The simplicial version of the $\star$ operation defined here is 
obviously functorial in the sense that \ 
$(e^{\ast}f)\star g= e^{\ast} (f\star g)$ \ for \ $e:F\C\ra F\A$ \ in
$\GG$, and \ $f\star (g^{\ast}h)=(f\star g)^{\ast}h$ \ for any $H$-map \ 
$h:\langle\K,m\rangle\ra\langle\LL,n\rangle$ \ between fibrant
$H$-simplicial sets which is strictly multiplicative (i.e., \ 
$n\circ(h\times h)= h\circ m:\K\times\K\ra\LL$).

However, Proposition \ref{pone} is still valid only in the homotopy
category, and this is in fact the obstruction to $\K$ being equivalent
to a loop space:
\end{remark}

%
%
\begin{thm}\label{tone}\stepcounter{subsection}
If $\K$ is an $H$-group in \ $\Sk$ \ such that \ 
$$
(\ast)\hsm f\star(g\star h)=(f^{\#}g)\star h \hsm
\forall f:F\A\ra F\B \ \mbox{and \ $g:F\B\ra F\C$ \ in $\GG$ 
and \ $h:\C\ra\K$},
$$
\noindent then $\K$ is $H$-homotopy equivalent to a simplicial group 
(and thus to a loop space); \ conversely, if \ $\K\in \GG$ \ 
(in particular, if \ $\K=G\LL$ \ for some \ $\LL\in\Sc$), \ then \ 
$(\ast)$ \ holds.
\end{thm}

\begin{proof}
Assume that $\K$ is an $H$-group in \ $\Sk$ \ satisfying \ ($\ast$). \ 
We shall need a simplicial variant of Stover's construction of 
resolutions by spheres (Def.\ \ref{drxs}), \ so as in \cite[\S 2]{StoV}, 
define a comonad \ $L:\GG\ra\GG$ \ by
\setcounter{equation}{\value{thm}}\stepcounter{subsection} 
\begin{equation}\label{eseven}
L\G = \ \coprod _{k=0}^{\infty } \ 
               \coprod _{\phi\in Hom_{\GG}(FS^{k},G)} \ F\bS{k}_{\phi} \ \ 
               \bigcup \ \ \coprod _{k=0}^{\infty } \ 
               \coprod _{\Phi\in Hom_{\GG} (F\be{k+1},G)} \ F\be{k+1}_{\Phi},
\end{equation}
\setcounter{thm}{\value{equation}}

\noindent where \ $F\be{k+1}_{\Phi}$, \ the $\GG$-disc indexed by \ 
$\Phi:F\be{k+1}\to \G$, \ is attached to \ $F\bS{k}_{\phi}$, \ 
the $\GG$-sphere indexed by \ $\phi=\Phi|_{F\partial \be{k+1}}$,  \ by
identifying \ $F\partial \be{k+1}$ \ with \ $F\bS{k}$ \ (see \S
\ref{ssg} above).  The coproduct here is just the (dimensionwise) free
product of groups; the counit \ $\eps:L\G\ra \G$ \ is 
``evaluation of indices'', and the comultiplication \ 
$\vartheta:L\G\hra L^{2}\G$ \ is as in \S \ref{dfpc}.

Now let
$$
\W = \ \bigvee _{k=1}^{\infty } \ 
               \bigvee _{f\in Hom_{\GG}(S^{k},K)} \ \bS{k}_{f} \ \ 
               \bigcup \ \ \bigvee_{k=1}^{\infty } \ 
               \bigvee_{F\in Hom_{\Ss} (e^{k+1},K)} \ \be{k+1}_{F}
$$

\noindent (the analogue for $\Ss$ \ of \ $L\G$, \ with the corresponding
identifications),  and let \ $z:\W\ra\K$ \ be the counit map. Then 
$z$ induces an epimorphism \ $z_{\ast}:\pis\W\epi\pis\K$ \ of \Pa s. \ 
($\K$ is a Kan complex, but $\W$ is not, so we understand \ $\pis\W$ \ 
to be the corresponding free \Pa\ $\cong \pis|\W|$ \ -- \ cf.\ \S \ref{dfpa}).

Likewise, we have an epimorphism of \Pa s \ 
$\tilde{\zeta}:\pis\Sigma\W\epi\Gs$, \ where \ 
$\Gs=\Omega^{-1}\pis\K$ \ is the delooping of \ $\pis\K$ \ -- \ or
equivalently, \ $\tilde{z}_{\ast}:\pis^{H}F\W\ra\pis^{H}\K$, \ induced
by \ $\tilde{z}=\bar{m}\circ Fz:F\W\ra\K$ \ (cf.\ \S \ref{dsv}).

Let \ $\M_{n}=L^{n}F\W$ \ for \ $n=0,1,\ldots$, \ with face and
degeneracy maps determined by the comonad structure maps \ $\eps$, 
$\vartheta$ \ -- \ except for \ $d_{n}:\M_{n}\ra\M_{n-1}$, \ 
defined \ $d_{n}=L^{n-1}\bar{d}$, \ where \ 
$\bar{d}:LF\W\ra F\W$, \ restricted to a summand \ 
$F\A_{\alpha}$ \ in \ $LF\W$ \ ($\A=\bS{k},\be{k+1}$), \ is an isomorphism 
onto \ $F\A_{\beta}\hra F\W$, \ where \ $\beta:\A\ra\K$ \ is the
composite \ $(\alpha\star z)\circ j_{A}$\vsm .

Because \ ($\ast$) \ holds exactly, we may verify that \ 
$\bar{d}\circ T\bar{d}=\bar{d}\circ T\eps:\M_{2}\ra \M_{0}$, \ so that \ 
$\Md$ \ is a simplicial object over $\GG$. \ Moreover, the augmented 
simplicial \Pa\ \ $\pis\Md\stackrel{\tilde{\zeta}}{\ra}\Gs$ \ is acyclic, 
by a variant of \cite[Prop.\ 2.6]{StoV}.
Thus the Quillen spectral sequence for \ $\Md$ \ (see \eqref{efour}) \ has \ 
$E^{2}_{s,t}=0$ \ for \ $s>0$, \ and \ $E^{2}_{0,\ast}\cong \Gs$, \ so it
collapses, and \ $\pis^{H}\diag\Md\cong\Gs=\pis^{H}\K$. \ 
Therefore, if we set \ $\LL=\diag\bar{W}\Md$ \ (which is isomorphic to \ 
$\bar{W}\diag\Ld$) \ we obtain a Kan complex $\LL$ such that \ 
$\K\simeq G\LL$ \ -- \ so \ $|\K|\simeq\Omega|\LL|$\vsm .

The converse is clear, since if \ $\K\in\GG$ \ then \ $j_{A}:\A\ra F\A$ \ 
induces a one-to-one correspondence between maps \ 
$f:\A\ra\K$ \ in \ $\Ss$ \ and homomorphisms \ $\varphi:F\A\ra\K$ \ 
in $\GG$, \ by the universal property of $F$.
\end{proof}

%
%
\section{Rectifying simplicial spaces}
\label{crss}

Theorem \ref{tone} suggests a way to determine whether an $H$-group $\X$ 
is equivalent to a loop space. Note that in fact we need only verify 
that \ \ref{tone}($\ast$) \ holds for \ $\A$, $\B$, and $\C$ \ in \ $\Ss$ \ 
which are homotopy equivalent to wedges of spheres. We now suggest a universal
collection of examples which may be used for this purpose, organized into one
(very large!) simplicial diagram. First, some definitions:

\begin{defn}\label{dsuh}\stepcounter{subsection}
A \textit{simplicial space up-to-homotopy} is a diagram \ $\hWd$ \ 
over \ $\Ts$ \ consisting of a sequence of spaces \ $\W_{0},\W_{1},\ldots$, \ 
together with face and degeneracies maps \ 
$d_{i}:\W_{n}\ra \W_{n-1}$ \ and \ $s_{j}:\W_{n}\ra \W_{n+1}$ \ 
($0\leq i,j\leq n$) \ -- \ satisfying the simplicial identities only 
\textit{up to homotopy}. 

Note that such a diagram constitutes an ordinary simplicial object over \ 
$ho\Ts$, \ so we can apply the functor \ $\pis:\Ts\ra\PAlg$ \ to \ 
$\hWd$ \ to obtain an (honest) simplicial \Pa\ \ 
$\pis(\hWd)\in s\PAlg$. \ Similarly for a \pss\ up-to-homotopy \ 
$\hWud$. 
\end{defn}

\begin{remark}\label{rcum}\stepcounter{subsection}
Note that diagrams denoted by \ $\Wd$, \ $\hWd$, \ $\Wud$, \ and  \ $\hWud$ \ 
each consist of a sequence of spaces \ $\W_{0},\W_{1},\ldots$; \ they differ 
in the maps with which they are equipped, and whether the identities which 
these maps are required to satisfy must hold in \ $\Ts$ \ or only in \ $ho\Ts$.
\end{remark}

\begin{defn}\label{drec}\stepcounter{subsection}
An (ordinary) simplicial space \ $\Vd\in s\Ts$ \ is called
a \textit{rectification} of a simplicial space up-to-homotopy \ $\hWd$ \ 
if \ $\V_{n}\simeq \W_{n}$ \ for each \ $n\geq 0$, \ and the face and 
degeneracy maps of \ $\Vd$ \ are homotopic to the corresponding maps of \ 
$\hWd$ \ (see \cite[\S 2.2]{DKSmH}, e.g., for a more precise definition).
For our purposes all we require is that \ $\pis\Vd$ \ be isomorphic (as a 
simplicial \Pa) to \ $\pis(\hWd)$. \ Similarly for rectification of 
\pss s, and \ ($\Delta$-)simplicial objects in \ $ho\Sk$ \ or \ $ho\GG$.
\end{defn}

\subsection{a \pss\ up-to-homotopy}
\label{spss}\stepcounter{thm}
Given an $H$-group $\X$, \ we wish to determine whether it is a loop space,
up to homotopy. We start by choosing some free simplicial \Pa\ \ $\Asd$ \ 
resolving \ $\Gs=\Omega^{-1}\pis\X$. \ 
By Remark \ref{rfpa}, the free simplicial \Pa\ \ $\Asd$ \ corresponds to
a simplicial object over the homotopy category, unique up to isomorphism 
(in \ $ho\Ts$), \ with each space homotopy equivalent to a wedge of spheres.
Therefore, it may be represented by a simplicial space up-to-homotopy \ 
$\hWd$, \ with \ $\pis(\hWd)\cong\Asd$ \ (\S \ref{dsuh}). We denote its 
underlying \pss\ up-to-homotopy by \ $\hWud$.

In light of Theorem \ref{tone} it would perhaps 
be more natural to consider the corresponding simplicial object up-to-homotopy 
over $\GG$, or $\Ss$, but given the equivalence of homotopy 
categories \ $ho\Ts\cong ho\GG\cong ho\Ss$, \ we prefer to work in the 
more familiar topological category.

Now \ $\hWd$ \ may be rectified if and only if it can be made 
$\infty$-\textit{homotopy commutative} \ -- \ that is, if and only if
one can find a sequence of homotopies for the simplicial identities
among the face and degeneracy maps, and then homotopies between these,
and so on (cf.\ \cite[Corollary 4.21 \& Theorem 4.49]{BVoHI}).
An obstruction theory for this was described in \cite{BlaHH}, and we 
briefly recall the main ideas here, mainly because we wish to present a 
technical simplification which eliminates the need for \cite[\S 6]{BlaHH}:
as we shall see below, it suffices to rectify \ $\hWud$; \ so we describe
an obstruction theory for the rectification of \pss s up-to-homotopy.
For this, we need some definitions from \cite[\S 5]{BlaHH}:

\begin{defn}\label{dperm}\stepcounter{subsection}
The $k$-dimensional \textit{permutohedron} \ $P_{k}$ \ is defined to be
the convex hull in \ $\R^{k+1}$ \ of the \ $(k+1)!$ \ points \ 
$(\sigma(1),\sigma(2),\ldots,\sigma(k+1))\in\R^{k+1}$, \ indexed by 
permutations \ $\sigma\in\Sigma_{k+1}$ \ (cf.\ \cite[0.10]{ZiegL}). Its
boundary is denoted by \ $\partial P_{k}$.

For \ $n\geq 0$ \ and any morphism \ 
$\delta:\bn$+$\bo\ra\bn$-$\bk$ \ in \ $\dDo$ \ 
(see \S \ref{dso} above), we may label the vertices of \ $P_{k}$ \ by all
possible ways of writing $\delta$ as a composite of face maps 
(cf.\ \cite[Lemma 4.7]{BlaHH}), \ and one can similarly interpret the faces 
of \ $P_{k}$. \ We shall write \ $P_{k}(\delta)$ \ for \ $P_{k}$ \ so labelled
(thought of as an abstract combinatorial polyhedron).
\end{defn}

\begin{defn}\label{dcc}\stepcounter{subsection}
Let \ $\hWud$ \ be a \pss\ up-to-homotopy, and \ 
$\delta:\bn$+$\bo\ra\bn$-$\bk$ \ some morphism in\ $\dDo$. \ We denote by \ 
$C(\delta)$ \ the collection of all proper factors of $\delta$ \ -- \ 
that is, \ 
$\gamma\in C(\delta)\Leftrightarrow \gamma'\circ\gamma\circ\gamma''=\delta$ \ 
and \ $\gamma',\gamma''$ \ are not both \ $id$.

A \textit{compatible collection for \ $C(\delta)$ \ and} \ $\hWud$ \ 
is a set \ $\{g^{\gamma}\}_{\gamma\in C(\delta)}$ \ of maps \
$g^{\gamma}:P_{n-k}\ltimes\W_{n}\ra\W_{k-1}$, \ one for each \
$\gamma=[(i_{k},\ldots,i_{n})]\in C(\delta)$, \ such that for any partition \ 
$\langle~i_{k},\ldots i_{\ell_{1}}~|~i_{\ell_{1}+1},\ldots i_{\ell_{2}}~|~
\ldots | i_{\ell_{r-1}+1},\ldots i_{n}~\rangle$ \ of \ 
$i_{k},\dotsc, i_{n}$ \ into $r$ blocks \ (where \ 
$\gamma=d_{i_{n}}\circ\dotsb\circ d_{i_{k}}$ \ in \ $\dDo$), \ with \ 
$\gamma_{1}=,d_{i_{\ell_{1}}}\circ\dotsc d_{i_{k}},\dotsc,
\gamma_{r}=d_{i_{n}}\circ\dotsb\circ d_{i_{\ell_{r-1}+1}}$, \ and we set \ 
$P\DEF P_{\ell_{1}-k-1}(\gamma_{1})\times P_{\ell_{2}-\ell_{1}}(\gamma_{2})
\dotsb P_{n-\ell_{r-1}}(\gamma_{r})$, \ 
then we require that \ $g^{\gamma}|_{P\ltimes Y{n}}$ \ be the 
composite of the corresponding maps \ $g^{\gamma_{i}}$, \ in the obvious sense.
We further require that if \ $\gamma=[i_{j}]$, \ then \
$g^{\gamma}$ \ must be in the prescribed homotopy class of \
$[d_{i_{j}}]\in[\W_{j+1},\ \W_{j}]$. 
\end{defn}

We shall be interested in such compatible collections only up to the obvious
homotopy relation. Note that for any \ $\delta:\bn$+$\bo\ra\bn$-$\bk$ \ in \ 
$\dDo$, \ any compatible collection \ $\{g^{\gamma}\}_{\gamma\in C(\delta)}$ \ 
induces a map \ $f=f^{\delta}:\partial P_{k}\ltimes\W_{n+1}\ra\W_{n-k}$, \ 
and compatibly homotopic collections induce homotopic maps.

\begin{defn}\label{dhho}\stepcounter{subsection}
Given \ $\hWud$ \ as in \S \ref{spss}, for each \ $k\geq 2$ \ and \ 
$\delta:\bn$+$\bo\ra\bn$-$\bk\in\dDo$,\ the $k$-{\em th order homotopy 
operation\/} (associated to \ $\hWud$ \ and $\delta$) \ is 
a subset \ $\lrc{\delta}$ \ of the track group \ 
$[\Sigma^{k-1}\W_{n+1},\W_{n-k}]$, \ defined as follows:

Let \ $S\subseteq [\partial P_{k}\ltimes\W_{n+1},~\W_{n-k}]$ \ 
be the set of homotopy classes of maps \ 
$f=f^{\delta}:\partial P_{k}(\delta)\ltimes\W_{n+1}\ra\W_{n-k}$ \ which are
induced as above by some compatible collection \ 
$\{g^{\gamma}\}_{\gamma\in C(\delta)}$. \ Choose a splitting \ 
$$
\partial P_{k}(\delta)\ltimes \W_{n+1}\cong 
\bS{k-1}\ltimes\W_{n+1}\simeq \bS{k-1}\wedge \W_{n+1}\vee \W_{n+1},
$$
and let \ $\lrc{\delta}\subseteq [\Sigma^{k-1}\W_{n+1},\ \W_{n-k}]$ \ 
be the image under the resulting projection of the subset \ 
$S\subseteq [\partial P_{k}\ltimes\W_{n+1},\W_{n-k}]$.
\end{defn}

\subsection{coherent vanishing}
\label{scv}\stepcounter{thm}

It is clearly a necessary condition in order for the subset \ 
$\lrc{\delta}$ \ to be non-empty that all the lower order operations 
(for \ $\gamma\in C(\delta)$) \text{vanish} \ -- \ i.e., 
contain the null class; a sufficent condition is that they do so 
\textit{coherently\/}, in the sense of \cite[\S 5.7]{BlaHH}. 
Again one may define a collection of higher homotopy operations in
various track groups \ $[\Sigma\W_{m},\W_{m-\ell}]$, \ 
whose vanishing guarantees the coherence of given collection of maps
(see \cite[\S 5.9]{BlaHH}). One then has
%
%
\begin{prop}[see Theorem 6.12 of \cite{BlaHH}]
\label{pthree}\stepcounter{subsection}
Given a \pss\ up-to-\-ho\-mo\-to\-py \ $\hWud$, \ it may be rectified to
a strict \pss\ \ $\Vud$ \ if and only if all the sequence of higher 
homotopy operations defined above vanish coherently.
\end{prop}

\subsection{adding degeneracies}
\label{sad}\stepcounter{thm}

Now assume given a \pss\ \ $\Wud\in\del\Ts$, \ to which we wish 
to add degeneracies, in order to obtain an full simplicial space \ $\Wd$.
Note that because the original \pss\ up-to-homotopy \ $\hWud$ \ of 
\S \ref{spss} above was obtained from the simplicial \Pa \ $\Asd$, \ in the case of interest to us \ $\Wud$ \ is already equipped with degeneracy maps \ -- \ 
but these satisfy the simplicial identities only up to homotopy!

In this situation, a similar obstruction theory was defined in 
\cite[\S 6]{BlaHH} for rectifying the degeneracies; but it was conjectured 
there that this theory is actually unnecessary (\cite[Conj.\ 6.9]{BlaHH}).
We now show this is in fact correct. 

\begin{defn}\label{dmat}\stepcounter{subsection}
Given a \pss\ $\Vud$, \ its $n$-th \textit{matching space} \ 
$M_{n}\Vud$ \ is defined to be the limit \ 
$$
M_{n}\Vud\DEF \{(x_{0},\dotsc,x_{n})\in(\V_{n-1})^{n+1}~\lvert\ 
d_{i}x_{j}=d_{j-1}x_{i} \text{\ \ for all\ }0\leq i<j\leq n\}.
$$
The map \ $\delta_{n}:\V_{n}\ra M_{n}\Vud$ \ is defined \ 
$\delta_{n}(x)=(d_{0}x,\dotsc,d_{n}x)$. \ 
(See \cite[XVII, 87.17]{PHirL}, and compare \cite[X,\S 4.5]{BKaH}).
\end{defn}

\begin{defn}\label{drf}\stepcounter{subsection}
A \pss\ $\Vud$ \ is called \textit{Kan} if for each \ $n\geq 1$ \ the
map \ $\delta_{n}:\V_{n}\ra M_{n}\Vud$ \ is a fibration. (See 
\cite[XVII, 88.2]{PHirL} or \cite[XII, \S 54]{DHKanM}, 
where this is called a Reedy fibrant object).
\end{defn}

%
%
\begin{lemma}\label{ltwo}\stepcounter{subsection}
For any \pss\ \ $\Xud$, \ there is a Kan \pss\ \ $\Vud$ \ and a map of
\pss s \ $f_{\bullet}:\Xud\ra\Vud$ \ such that each \ $f_{n}$ \ is a
homotopy equivalence.
\end{lemma}

\begin{proof}
This follows from the existence of the so-called Reedy model category
structure on \ $\del\Ts$ \ (see \cite[XVII, Thm.\ 88.3]{PHirL}); \ 
$\Vud$ \ may be constructed directly from \ $\Xud$ \ by successively
changing the maps \ $\delta_{n}$ \ into fibrations as in
\cite[I, (7.30)]{GWhE}.
\end{proof}

Thus the proof of \cite[Conj.\ 6.9]{BlaHH} follows from:

%
%
\begin{prop}[compare Theorem 5.7 of \cite{RSaD1}]
\label{pfour}\stepcounter{subsection}
If \ $\Vud$ \ is a Kan \pss\ which rectifies the \pss\ up-to-homotopy \ 
$\hWud$ \ of \S \ref{spss}, then one can define degeneracy maps on \ 
$\Vud$ \ making it into a full simplicial space.
\end{prop}

\begin{proof}
Using the singular functor \ $S:\Ts\ra\Sk$ \ we may work with simplicial sets,
rather than topological spaces; the maps \ 
$\delta_{n}:\V_{n}\ra M_{n}\Vud$ \ are now assumed to be Kan fibrations 
in \ $\Ss$.

By induction on \ $n\geq 0$ \ we assume that degeneracy maps \ 
$s_{j}:\V_{k}\ra\V_{k+1}$ \ have been chosen for all \ 
$0\leq j\leq k<n$, \ satisfying all relevant simplicial identities. 

Let \ $\Sigma$ \ denote the subcategory of \ $\Do$ \ with the same objects 
as \ $\Do$, \ but only the degeneracies as morphisms. For each \ $k\geq 0$, \ 
let \ $\Sigma/\bk$ \ denote the ``over category'' of $\bk$. \ By assumption \ 
$\Vud$, \ together with the existing degeneracies, defines a functor \ 
$V:\Sigma/\bn$-$\bo\ra \Ss$. \ Denote its colimit by \ 
$L_{n}$; \ this may be thought of as  the sub-simplicial set 
of \ $\V_{n}$ \ consisting of the degenerate simplices; i.e., \ 
$\bigcup_{j=0}^{n-1} Im(s_{j})$. \ We also have an associated ``free'' 
functor \ $\Sigma/\bm\ra\Ss$ \ for each \ $m\geq 0$ \ (with the same 
values on objects as $V$, but all morphisms trivial). \ Its colimit, denoted 
by \ $D_{m}$, \ is the coproduct (i.e., wedge) over \ $Obj(\Sigma/\bm)$ \ of 
the spaces \ $\V_{k}$ \ ($0\leq k<m$), \ indexed by all possible iterated 
degeneracies \ $s_{i_{1}}\circ\dotsb\circ s_{i_{m-k}}:\bk\ra\bm$. \ 
This comes equipped with structure maps \ $e_{j}^{m-1}:\V_{m-1}\ra D_{m}$ \ 
($0\leq j<m$). \ 
See \cite[p.\ 95]{MayS} or \cite[\S 4.5.1]{BlaH} for an explicit description. 
For example, \ $L_{0}=D_{0}=\ast$, \ $L_{1}=D_{1}\cong \V_{0}$, \ but \ 
$D_{2}=(\V_{1})_{s_{0}}\vee(\V_{0})_{s_{1}s_{0}}\vee(\V_{1})_{s_{1}}$, \ 
while \ 
$L_{2}=(\V_{1})_{s_{0}}\amalg_{(\V_{0})_{s_{1}s_{0}}} (\V_{1})_{s_{1}}$ \ 
(the pushout).

In fact, if we define \ $\Yd\in s\Ss$ \ by \ $\Y_{n} \DEF \V_{n}\vee D_{n}$, \ 
with the obvious degeneracies (defined by the structure maps \ $e_{j}^{m}$) \ 
and face maps (induced from those of \ $\Vud$ \ via the simplicial 
identities), then \ $F(\Vud)\DEF \Yd$ \ defines a functor \ 
$F:\del\Ss\ra s\Ss$ \ which is left adjoint to the forgetful functor \ 
$U:s\Ss\ra\del\Ss$, \ and there is a natural inclusion \ 
$\iota:\Vud\ra UF\Vud$. \ 

Note that the degeneracies up-to-homotopy \ $s_{j}':\V_{n}\ra\V_{n}$, \ 
which exist because \ $\Vud$ \ rectifies \ $U(\hWd)$ \ (where \ 
$\pis(\hWd)\cong\Asd$ \ as simplicial \Pa s), \ 
define a map \ $\sigma_{n+1}':D_{n+1}\ra V_{n+1}$.

Our object is to define inductively a retraction \ $\sigma:UF\Vud\ra\Vud$ \ 
in \ $\del\Ss$, \ starting with \ $\sigma_{0}=id_{\V_{0}}$, \ such that \ 
$\sigma_{n}\sim \sigma_{n}'$ \ for all $n$. \ This map $\sigma$ must 
commute with the degeneracies defined so far: that is, at the $n$-th stage 
we must choose a map \ $\sigma_{n+1}:D_{n+1}\ra\V_{n+1}$ \ homotopic to \ 
$\sigma_{n+1}'$, \ and then define \ $s_{j}^{n}:\V_{n}\ra\V_{n+1}$ \ by \
%
\setcounter{equation}{\value{thm}}\stepcounter{subsection} 
\begin{equation}\label{eeleven}
s_{j}^{n}\DEF\sigma_{n+1}\circ e_{j}^{n}\circ\iota_{n}.
\end{equation}
\setcounter{thm}{\value{equation}}

Moreover, together with the face maps of \ $\Vud$, \ the degeneracies 
chosen so far determine a map \ $\rho_{n}:D_{n+1}\Vud\ra M_{n+1}\Vud$, \ 
by the universal properties of the limit and colimit; the 
simplicial identities in \ $\hWd$ \ guarantee that \ 
$\delta_{n+1}\circ \sigma_{n+1}' \sim \rho_{n+1}$. \ Note than in order
the simplicial identities \ 
$d_{i} s_{j}=s_{j-1} d_{i}$ \ (for \ $i<j$),
$d_{j} s_{j}= d_{j+1} s_{j}=id$, \ and \ $d_{i} s_{j}=s_{j} d_{i-1}$ \ 
(for \ $i>j+1$) \ to be satisfied, it suffices that
%
\setcounter{equation}{\value{thm}}\stepcounter{subsection} 
\begin{equation}\label{enine}
\delta_{n+1}\circ \sigma_{n+1}=\rho_{n+1},
\end{equation}
\setcounter{thm}{\value{equation}}
that is, \ $\sigma_{n+1}$ \ must be a lift for the given map \ $\rho_{n+1}$. \ 
On the other hand, in order that \ 
$s_{j} s_{i}= s_{i+1} s_{j}$ \ hold for all \ $j\leq i$, \ it suffices to have
%
\setcounter{equation}{\value{thm}}\stepcounter{subsection} 
\begin{equation}\label{eten}
\sigma_{n+1}\circ e^{n}_{j} = 
\sigma_{n+1}\circ e^{n}_{j} \circ \iota_{n}\circ \sigma_{n}
\hsp\text{for all\hsm }0\leq j\leq n,
\end{equation}
\setcounter{thm}{\value{equation}}
\noindent (where \ $\iota_{n}:\V_{n}\hra (UF\Vud)_{n}$ \ is the 
inclusion)\vsm .

Now \ $D_{n+1}$ \ has a wedge summand \ $\bar{D}_{n+1}$ \ such that \ 
$D_{n+1}=\bar{D}_{n+1}\vee\bigvee_{j=0}^{n}(\V_{n})_{s_{j}}$, \ 
and \ $\sigma_{n+1}':D_{n+1}\ra \V_{n+1}$ \ thus defines a map \ 
$\bar{\sigma}_{n+1}'=
\sigma_{n+1}'\rest{\bar{D}_{n+1}}:\bar{D}_{n+1}\ra \V_{n+1}$. \ 
Since \ $\delta_{n+1}$ \ is a fibration, one may use the homotopy lifting 
property to obtain a map \ $\bar{\sigma}_{n+1}\sim\bar{\sigma}'_{n+1}$ \ 
such that \ 
$\delta_{n+1}\circ \bar{\sigma}_{n+1}=\rho_{n+1}\rest{\bar{D}_{n+1}}$. \ 

Note that \ $L_{n}=Im(\sigma_{n})$, \ by induction, so \eqref{eten} for \ 
$n-1$ \ and \eqref{eeleven} imply that, for each \ $0\leq j\leq n$, \ 
the map \ $\bar{\sigma}_{n+1}\circ e_{j}^{n}:D_{n}\ra\V_{n+1}$ \ 
induces a map \ $g_{j}:L_{n}\ra\V_{n+1}$.

Because \ $\Asd$ \ is a free simplicial \Pa\ (Def.\ \ref{dfsr}) and \ 
$\pis(\hWd)\cong\Asd$, \ Lemma \ref{lthree} guarantees that there is a \ 
$\Z_{n}\in\Ss$, \ weakly equivalent to a wedge of spheres, and a map \ 
$f_{n}:\Z_{n}\ra\V_{n}$ \ which, 
together with the inclusion \ $h_{n}:L_{n}\hra\V_{n}$
induces a weak equivalence of simplicial sets \ 
$(h_{n},f_{n}):L_{n}\vee\Z_{n}\xrightarrow{\simeq}\V_{n}$. \ Since \ 
$h_{n}$ \ is a cofibration, \ using a minimal complex for \  
$\Z_{n}$ \ (see \cite[\S 9]{MayS}) we may assume that \ $(h_{n},f_{n})$ \ 
is a trivial cofibration in \ $\Ss$ \ (cf.\ \cite[II, 3.14]{QuH}). 
Again the fact that \ $\delta_{n+1}$ \ is a fibration implies that there exists
a lifting \ $\alpha_{j}:\Z_{n}\ra\V_{n+1}$ \ for \ 
$\rho_{n+1}\rest{(\V_{n})_{s_{j}}}\circ f_{n}:\Z_{n}\ra M_{n+1}\Vud$. \ 
Thus the left lifting property (cf.\ \cite[I, 5.1]{QuH}) for the solid
commutative square

%
%
%
\begin{picture}(170,110)(-100,-15)
\put(20,60){$L_{n}\vee\Z_{n}$}
\put(65,65){\vector(1,0){82}}
\put(82,71){$(g_{j},\alpha_{j})$}
\put(150,60){$\V_{n+1}$}
 
\put(40,57){\vector(0,-1){46}}
\put(0,32){$(h_{n},f_{n})$}
\put(-50,32){triv cof.}
\put(30,0){$(\V_{n})_{s_{j}}$}
 
\put(155,57){\vector(0,-1){46}}
\put(158,32){$\delta_{n+1}$\ \ fib.}
\put(145,0){$M_{n+1}\Vud$}
\put(65,5){\vector(1,0){72}}
\put(80,13){$\rho_{n+1}\rest{(\V_{n})_{s_{j}}}$}
 
\multiput(57,12)(6,3){15}{\circle*{.3}}
\put(145,57){\vector(2,1){2}}
\put(60,40){$(\sigma_{n+1})_{s_{j}}$}
 
\end{picture}

\noindent guarantees the existence of a dotted lifting \ 
$(\sigma_{n+1})_{s_{j}}:(\V_{n})_{s_{j}}\ra \V_{n+1}$ \ for \ 
$\rho_{n+1}\rest{(\V_{n})_{s_{j}}}$, \ 
and these liftings, for various $j$, together with \ $\bar{\sigma}_{n+1}$, \ 
define \ $\sigma_{n+1}:D_{n+1}\ra\V_{n+1}$ \ satisfying 
\eqref{eten} (and of course \eqref{enine}), as required. \ \eqref{eeleven} 
then defines \ $s_{j}:\V_{n}\ra\V_{n+1}$ \ for all \ $0\leq j\leq n$, \ 
completing the induction.
\end{proof}

%
%
\begin{cor}\label{cthree}\stepcounter{subsection}
If \ $\Vud$ \ is a \pss\ such that \ $\pis\Vud$ \ is a free simplicial \Pa\ 
(Def.\ \ref{dfsr}), there is a spectral sequence with \ 
$E^{2}_{s,t}=\pi_{s}(\pi_{t}\Vud) \Rightarrow \pi_{s+t}\mreal{\Vud}$.
\end{cor}

\begin{proof}
See \S \ref{sreal} and \eqref{ethree}, noting that the definition of the 
homotopy groups of a simplicial group is also valid for a $\del$-simplicial
group (see \cite[\S 17]{MayS}), and that in the proof of Proposition 
\ref{pfour} we did not use the fact that \ $\Asd$ \ was a resolution of \ 
$\Gs$.
\end{proof}

If the higher homotopy operations described in \S \ref{dhho} vanish 
coherently, then the \pss\ up-to-homotopy \ $\hWud$ \ of \S \ref{spss} 
may be rectified to a strict \pss\ \ $\Wud$, \ which may in 
turn be replaced by a Kan \pss\ \ $\Vud$ \ using Lemma \ref{ltwo}, 
with \ $\pis\Vud\cong\Asd$. The spectral sequence of Corollary \ref{cthree}
then implies that \ $\Y\DEF\mreal{\Vud}$ \ satisfies: \ 
$\pis\Y\cong \Gs=\Omega^{-1}\pis\X$. \ We have thus realized the algebraic 
delooping of \ $\pis\X$ \ by a space $\Y$.

\begin{remark}\label{rot}\stepcounter{subsection}
As in any obstruction theory, if one of the homotopy operations on
question does not vanish (or if there is a non-vanishing obstruction 
to coherence, as in \S \ref{scv}), one must backtrack, changing choices made
in previous stages. On the face of it, if all such choices show that the
\pss\ up-to-homotopy \ $\hWud$ \ cannot be rectified, we must then try 
other choices for the resolution \ $\Asd\ra\Gs$. \ However, 
we conjecture that in fact if \textit{one} free simplicial \Pa \ 
resolution of \ $\Gs=\Omega^{-1}\pis\X$ \ is realizable, then \textit{any} 
resolution is realizable (so that that any obstruction to rectifying \ 
$\hWd$ \ shows that $\X$ is not a loop space).
\end{remark}

\subsection{realizing \Pa\ morphisms}
\label{srpm}\stepcounter{thm}

It remains to ascertain that the space $\Y$ which realizes \ $\Gs$ \ is in 
fact a delooping of $\X$. In other words, we have an abstract 
\Pa\ isomorphism \ $\phi:\pis\Omega\Y \xrightarrow{\cong}\pis\X$ \ 
(cf.\ Corollary \ref{cone}), which we wish to realize by a map of spaces \ 
$f:\Omega\Y \ra \X$. \ Now, there is an obstruction theory for the realization 
of \Pa\ morphisms, simpler than but similar in spirit to that described above, 
which we briefly recapitulate. For the details, see \cite[\S 7]{BlaHH}, as 
simplified in \cite[\S 4.9]{BlaHO} (and see \cite[\S 4]{BlaHR} for an 
algebraic version).

We start with some \pas\ resolution of \ $\Omega\Y$ \ by wedges of 
spheres \ --  \ i.e., an augmented \pss\ \ 
$\Vud\xrightarrow{\varepsilon}\Omega\Y$ \ such that \ $\pis\Vud$ \ is a 
\pas\ \Pa\ resolution of \ $\pis\Omega\Y\cong\pis\X$, and each \ $\V_{n}$ \ 
is homotopy equivalent to a wedge of spheres (see \cite[\S 1]{StoV}). \ 
The spectral sequence of Corollary \ref{cthree} then implies that \ 
$\real{\Vud}\simeq \Omega\Y$. 

By Fact \ref{rfpa}, we can realize \ $\varepsilon:\pis\V_{0}\ra\pis\X$ \ by \ 
a map \ $e_{0}:\V_{0}\ra\X$, \ and then define \ $e_{n}:\V_{n}\ra\X$ \ by \ 
$e_{n}\DEF e_{n-1}\circ d_{n}$ \ for \ $n>0$. \ By the simplicial 
identities for \ $\pis\Vud\ra\pis\Omega\Y$, \ we know \ 
$\pis(e_{n})= \pis(e_{n-1})\circ d_{i}$, \  so that \ 
$e_{n}\sim e_{n-1}\circ d_{i}$ \ for all \ $0\leq i\leq n$. \ 
If we can make this hold on the nose, rather than just up to homotopy, 
then \ $\Vud\xrightarrow{e_{0}}\X$ \ is also a (strict) augmented \pss, so
the spectral sequence of Corollary \ref{cthree} now implies that \ 
$\real{\Vud}\simeq \X$, \ and thus  \ $\Omega\Y\simeq\X$. \ 
This is where the appropriate higher homotopy operations (defined as follows) 
come in:

\begin{defn}\label{dns}\stepcounter{subsection}
Let \ $D[n]$ \ denote the standard simplicial $n$-simplex, \ 
together with an indexing of its non-degenerate $k$-dimensional faces \ 
$D[k]^{(\gamma)}$ \ by the composite face maps \ 
$\gamma=d_{i_{n-k}}\circ\ldots\circ d_{i_{n}}:
\mathbf{n}\ra\mathbf{k}-\mathbf{1}$ \ in \ $\Do$ \ Its \ $(n-1)$-skeleton, 
which is a simplicial \ $(n-1)$-sphere, is denoted by \ $\partial D[n]$. \ 
We choose once and for all a fixed isomorphism \ 
$\varphi^{(\gamma)}:D[k]^{(\gamma)}\ra D[k]$ \ for each face \ 
$D[k]^{(\gamma)}$ \ of \ $D[n]$.
\end{defn}

\begin{defn}\label{dcs}\stepcounter{subsection}
Given \ $\Vud$ \ as above, for each \ $n\in\N$ \ we define 
a \ $\partial D[n]$-\textit{compatible sequence} to be a sequence of maps \ 
$\{ h_{k}:D[k]\ltimes \V_{k}\ra \X \}_{k=0}^{n-1}$, \ 
such that \ $h_{0}=e_{0}$, \ and for any iterated
face maps \ $\delta=d_{i_{j+1}}\circ\dotsb\circ d_{i_{n}}$ \ and \ 
$\gamma=d_{i_{j}}\circ\delta$ \ ($0\leq j<n$) \ we have \ 
$h_{j}\circ (id\ltimes d_{i_{j}}) = h_{j+1}\circ
(\iota^{\gamma}_{\delta}\ltimes id)$ \ on \ $D[j]\ltimes\V_{j+1}$, \ 
where \ $\iota^{\gamma}_{\delta}\DEF
\varphi^{\delta}\circ\iota\circ(\varphi^{\gamma})^{-1}$, \ and \ 
$\iota:D[j]^{(\gamma)}\ra D[j+1]^{(\delta)}$ \ is the inclusion.
\end{defn}

Any such a $\partial D[n]$-compatible sequence \ 
$\{ h_{k}:D[k]\ltimes\V_{k}\ra \X\}_{k=0}^{n-1}$ \ induces a map \ 
$\bar{h}:\partial D[n]\ltimes\V_{n}\ra\X$, \ defined on the ``faces'' \ 
$D[n-1]^{(d_{i})}\ltimes\V_{n}$ \ by \ 
$\bar{h}\rest{D[n-1]^{(d_{i})}\ltimes\V_{n}}=h_{n-1}\circ(id\ltimes d_{i})$, 
and for each \ $n\geq 2$, \ the $n$-{\em th order homotopy operation\/} 
(associated to \ $\Vud$) \ is a subset \ $\lrc{n}$ \ of the track group \ 
$[\Sigma^{n-1}\V_{n},\X]$ \ defined analogously to \S \ref{dhho}.

Again as in \S \ref{scv}, the coherent vanishing of all the operations \ 
$\{\lrc{n}\}_{n=2}^{\infty}$ \ is the necessary and sufficient condition 
for \ $\Omega\Y\simeq\X$. 

\begin{remark}\label{rshs}\stepcounter{subsection}
We observe that if we choose to work in the category \ $\Sk$, \ rather than \ 
$\Ts$, \ we replace \ $\Omega\Y$ \ by \ $G(S\Y)$, \ and we may then use the 
resolution \ $\Md\ra G(S\Y)$ \ of Theorem \ref{tone} (or rather, \ 
$\Md^{\Delta}\ra G(S\Y)$) \  instead of \ $\Vud\ra\Omega\Y$. \ 
The obstruction theory that arises is of course equivalent to the one we just 
sketched, but it fits directly into the philosophy of section \ref{csg}, 
since our theory implies that $\X$ is a loop space if and only if the 
$H$-group augmentation up-to-homotopy from the given bisimplicial group \ 
$\Md$ \ to \ $S\X$ \ can be rectified.
\end{remark}

We summarize the results of this section in
%
%
\begin{thm}\label{ttwo}\stepcounter{subsection}
An $H$-group $\X$ is $H$-equivalent to a loop space if and only if 
the collection of higher homotopy operations defined in \S \ref{dhho} and
\S \ref{srpm} above (taking value in homotopy groups of spheres or in \ 
$\pis\X$, \ respectively) vanish coherently. 
\end{thm}

\begin{remark}\label{ran}\stepcounter{subsection}
The obstruction theory described here is not applicable as such to the 
related question of the existence of \ $A_{n}$-structures on an $H$-space 
$\X$ (cf.\ \cite[\S 2]{StaH1}). An alternative approach to the loop-space 
question, more closely related to Stasheff's theory, but still expressable 
in terms of homotopy operations taking values in homotopy groups (rather 
than higher homotopies), will be described in a future paper, with a view 
to such an extension.
\end{remark}
%
%

\end{document}